\documentclass[11pt]{amsart}
\usepackage{amssymb}
\numberwithin{equation}{section}

\begin{document}

\title{Why Sam Northshield's recent one-line proof of the infinitude of primes is flawed?}

\author{O.A.S. Karamzadeh}

\address{Department of Mathematics, Shahid Chamran University of Ahvaz, Ahvaz, Iran}

\email{karamzadeh@ipm.ir}

 \maketitle

Recently a Note with the title "A one-line proof of the infinitude of primes" has appeared in Monthly's 2015 MathBit, see \cite{sam}.  A video of this proof has already gone viral on YouTube and some people are admiring this proof on the Internet, and there are also some who show disapproval and express doubts on the validity of this proof. Confer also \cite{yu}, to see a detailed proof of this one-line proof, which is still a flawed-proof, by the reasoning which follows, briefly. The author of \cite{yu}, also makes the comment that the one-line proof in \cite{sam} is clear and elegant, but to explain it to high-school students, it needs to be deciphered. I must remind the reader that, all that is needed, to be explained about this proof, in order to be possibly deciphered, is already done in a best possible way in the above video. This proof (together with its presentation in this video) might look fine to most of the high-school students and even to their teachers, for they might not have the needed experience to detect the existing  flaw in the proof (indeed, this is the main reason for writing this Note). Unfortunately, I must claim that this proof is certainly flawed and it is not decipherable in a positive direction.  Before, giving our detail reasoning, I must anticipate to bring to the reader's attention that the essence of the proof is what follows: in this one-line proof, simply an obvious contradiction is used to get another contradiction, which is all that is needed for this particular proof. Whereas we all know that no result (even an obtained contradiction) is valid if it is obtained by the assumption of a contradiction. We should bear in mind that all the flawed-proofs in the literature, which are created intentionally, are to amuse and to provoke the reader. These kind of proofs are usually published in certain magazines with a special column for these proofs. We should remind the reader that, the College Mathematics Journal, used to have a column with the title ``Fallacies, Flaws, and Flimflam, which was very useful to both the students and the teachers alike. Certainly as we all know, the American Mathematical Monthly does not have such a column devoted to these kind of proofs. Hence, this particular one-line proof in the Monthly is not to be considered as a flawed-proof or to be of any amusement for the reader. However, in my opinion, this proof is certainly flawed. Before going on and giving my reasons to refute the proof, I apologize in advance for repeating some unncessary details in my arguments.  I should also admit that this rather lengthy reasoning which seems to be, much ado about nothing, may not be needed.  But the reader might agree with me that these kind of discussions require scrupulous attentions to detail (especially, when the proof has already received much attention). Some readers might be led on by the "visual nature" of the proof and consider it as a proof without words (although, to comprehend the proof, even in that case, one has to repeat a kind of Euclid's proof in one's mind without any success, for completing the last part, i.e., getting the zero equality in the proof). But, if only one tries to write down a proof, a stumbling block might appear, with no way of getting around it. Not only this incorrect proof is read by many young students, it is also watched by millions of viewers who might share it with others to get the proof widespread. The fact that this might mislead many students and teachers (from mathematics education point of view), in schools throughout the world, has prevailed upon me to give some comments on this incorrect proof, which follows. 
Suppose that there are only finitely many primes and let P be their product and $F$ be the set of these primes. Then the one-line proof in \cite{sam} says: 

$0<\prod\limits_{p\in F}\sin\frac{\pi}{p}=\prod\limits_{p\in F}\sin \pi\frac {(1+2P)}{p}=0$.  No explanation is given.  Apparently, as with many of the Monthly's MathBits, the reader is left to supply the details. I believe some correct part of the the proof, on which, every reader is certainly agreed, may go like this: the left hand side product is manifestly positive and equals to the right hand side product.  But to see that $\prod\limits_{p\in F}\sin \pi\frac {(1+2P)}{p}=0$, we notice the latter product is zero if and only if one of the factors of this product is zero (note, this part is also manifest). Consequently, in order to get zero for the above product, at least  one of the factors of $\prod\limits_{p\in F}\sin \pi\frac {(1+2P)}{p}$ must be zero. But $\sin q\pi=0$, where $q$ is a rational number, if and only if $q$ is an integer (sorry again, to repeat trivialities).  This means that $\frac{1+2P}{p}$ must be an integer for some $p\in F$ (note, $p$ is a prime divisor of $P$). Consequently, we must somehow assume that $1+2P$ is divisible by a member, $p$ say, of $F$. This implies that $1+2P=kp$ for some $k$, i.e., $\sin\pi\frac {(1+2P)}{p}=\sin k\pi=0$. But at the same time we already know that $1+2Q$ cannot be divisible by any factor of $Q$, where $Q$ is any positive integer (note, the latter fact does not depend on whether the set of all primes is finite or infinite). In fact in some variants of the classical Euclid's proof, by making use of the fact that  $1+2P$ cannot be divisible by any factor of $P$, one gets the desired contradiction. Therefore, either we have to, give in, to the fact that $1+2P$ is never divisible by any factor of $P$, a fortiori, by any prime factor of $P$, i.e., by any member of $F$, in which case, it implies  that $\prod\limits_{p\in F}\sin \pi\frac {(1+2P)}{p}\ne 0$. In consequence, the one-line proof in [1] is incorrect. Or, ignoring this fact for a moment, and simply try to give a different proof, by some other reasoning, to show that the above product is zero. But the only possible way left to anyone to provide any arguments for such a proof (if possible?) is this: since we are assuming the set of all primes is $F$  and $1+2P$ is a positive integer, therefore we are forced, by the Fundamental Theorem of Arithmetic, to assume that $1+2P$ must be divisible by a member of $F$ (i.e., by a prime factor, $p$ say, of $P$). But this compulsory assumption immediately leads us to a contradiction (note, $1+2P$ is never divisible by a member of $F$), therefore we must pause and avoid deducing any new facts based on this contradictory assumption. In particular, we should avoid inferring that $\sin\pi \frac{1+2P}{p}=0$ (note, we should emphasize again the fact that $1+2P$ is never divisible by a prime factor of $P$ is independent of our assumption that the set of primes is finite, and therefore cannot be altered by our above temporary ignorance). We should also emphasize, in any proof, when we reach a contradiction, we must pause for a moment and correct the false assumption that causes this contradiction, before continuing the proof with this contradiction. This means that  the compulsory assumption, which has already caused a contradiction, leaves us with no other choice but to assume that the set $F$ must be infinite before doing anything else (i.e., correcting the false assumption before continuing the proof). Hence, the proof of the infinitude of primes is complete at this stage (i.e., in the middle of the one-line proof we are done). We should remind the reader what happens, in the previous argument
 (i.e., completion of the proof at the above stage), is not because of the one-line proof in [1].  In fact, we couldn't still get zero for the above product in the one-line proof. It is indeed the Euclid's proof, in that argument, which completes the proof halfway. We must repeat ,with emphasis, again that we cannot go on (like the only possible way left for the proof in [1]), with this compulsory assumption, which is an obvious contradiction. For otherwise, if we go on with the compulsory assumption (i.e., with a contradiction), which is clearly the essence of the proof in [1], it means that we couldn't avoid using this contradiction before finishing this proof. This is where the "{\bf principle of explosion}" comes in, which is the principle of classical logic that states anything follows from a contradiction. That is to say, once a contradiction has been assumed, any proposition (including their negations) can be deduced from it, i.e., one can even assume the above product is also nonzero. This means that if we ignore the above fact (i.e., the fact that $1+2P$ is never divisible by any prime factor of $P$), the one-line proof becomes flawed. To me this does not differ from the well-known, flawed trivial proof, of showing any two real numbers, $a, b$ say, are equal, by ignoring the fact that real numbers are not divisible by zero. In sum, there is no way for any reader to provide details for the one-line proof in [1], in order to become a flaw-less proof. Let us conclude this Note with the comments which follows.  Consider the equality $\prod\limits_{p\in F}\pi\frac {(1+2P)}{p}=\prod\limits_{p\in F}\pi(\frac{1}{p}+\frac{2P}{p})$. Manifestly for all $p\in F$, none of the rational numbers $\frac{1+2P}{p}=\frac{1}{p}+\frac{2P}{p}$ is an integer, for $\frac{2P}{p}$ is an integer. Consequently,  for all the reasons on earth no one can infer from the equality of the above two products that $1+2P$ is divisible by a prime divisor of $P$ (note, all one may infer, of course, not necessarily as a consequence of the above equality, is the fact that, $1+2P$ is divisible by a prime number which cannot be in $F$). But, why on earth  should anyone believe that if we insert the sine function, artificially, into the above products, i.e.,  $\prod\limits_{p\in F}\sin\pi\frac {(1+2P)}{p}=\prod\limits_{p\in F}\sin\pi(\frac{1}{p}+\frac{2P}{p})=\prod\limits_{p\in F}\sin\frac{\pi}{p}>0$, then $\frac {1+2P}{p}=\frac{1}{p}+\frac {2P}{p}$ becomes an integer for a prime divisor, $p$ say, of $P$, and hence $\prod\limits_{p\in F}\sin\pi\frac {(1+2P)}{p}=0$ (note, this is taken for granted after applying the sine function in the one-line proof, which is an obvious flaw). What property of the sine function is responsible that makes $1+2P$ to become
 divisible by a prime factor of $P$? Aren't we  after obtaining a contradiction, artificially? By what we have already observed, it should also be emphasized that, contrary to the comment in \cite{yu}, this one-line proof cannot be deciphered either. Finally it seems, a flaw occurs in the proof in [1], because in this one-line proof we are somehow trying to make the Euclid's classical proof concealed.

\end{document}